\documentclass[11pt, twoside, eqno]{article}
\usepackage{latexsym}
\usepackage{amssymb}
\usepackage{amsfonts}
\begin{document}
\begin{center}

\noindent {\bf \Large Comment on ``Filters in ordered 
$\Gamma$-semigroups"}\bigskip

\medskip

{\bf Niovi Kehayopulu}\bigskip

{\small Department of Mathematics,
University of Athens \\
15784 Panepistimiopolis, Athens, Greece \\
email: nkehayop@math.uoa.gr }
\end{center}
\date{ }

\bigskip

\noindent{\bf ABSTRACT.} This is about the paper in the title by 
Kostaq Hila in Rocky Mt. J. Math. 41, no. 1 (2011), 189--203 [1] for 
which corrections should be done.\bigskip

\bigskip

\bigskip

\noindent Throughout the paper in [1], $M$ stands for an ordered 
$\Gamma$-semigroup (: $po$-$\Gamma$-semigroup) [1, page 191].

If $A$ is a left (or right) ideal of a $po$-$\Gamma$-semigroup $M$, 
then $A$ is a subsemigroup of $M$. So the ideals of $M$ are 
subsemigroups of $M$ as well. In the proof of Theorem 2.5 in [1], it 
has been proved that $(a)_{\cal N}$ is a subsemigroup and at the same 
time an ideal of $N(a)$. Since the ideals are subsemigroups, the 
proof that $(a)_{\cal N}$ is a subsemigroup of $N(a)$ should be 
deleted from the proof of the theorem and keep only the fact that 
$(a)_{\cal N}$ is a nonempty subset of $N(a)$. It might be emphasized 
in the theorem that the $N(a)$ is a $po$-$\Gamma$-semigroup.\\
This is the Theorem 2.6 in [1]:\smallskip

\noindent{\bf Theorem 2.6.} {\it Let $a,b\in M$. Then $(a)_{\cal 
N}\preceq (b)_{\cal N}$ if and only if $N(b)\subseteq 
N(a)$.}\smallskip

\noindent In the proof of the ``$\Rightarrow$'' part of the theorem 
``it is clear that $(a)_{\cal N}\preceq (b)_{\cal N}$ implies 
$(a)_{\cal N}\gamma (b)_{\cal N}\subseteq (a)_{\cal N}$ $\forall\, 
\gamma\in\Gamma$'' is written. The correct is that $(a)_{\cal N} 
\preceq (b)_{\cal N}$ implies $(a)_{\cal N}\gamma (b)_{\cal 
N}=(a)_{\cal N}$ $\forall\,\gamma\in\Gamma$. For the class $(a)_{\cal 
N}$ instead of using the (undefined concept) ``$(a)_{\cal N}$ is a 
semilattice congruence class'' the phrase ``$\cal N$ is a semilattice 
congruence on $M$'' should be used.

In the proof of the ``$\Leftarrow$'' part of the theorem ``we only 
need to prove that $(a)_{\cal N}\gamma (b)_{\cal N}\subseteq 
(a)_{\cal N}$, $(b)_{\cal N}\gamma (a)_{\cal N}\subseteq (a)_{\cal 
N}$ $\forall\,\gamma\in\Gamma$" is written. As $\cal N$ is a 
semilattice congruence on $M$, $(a)_{\cal N}\gamma (b)_{\cal 
N}=(b)_{\cal N}\gamma (a)_{\cal N}$. So the author has only proved 
that $(a)_{\cal N}\gamma (b)_{\cal N}\subseteq (a)_{\cal N}$ 
$\forall\,\gamma\in\Gamma$ and to prove the theorem, he had to prove 
that for all $\gamma\in\Gamma$, the inclusion $(a)_{\cal N}\subseteq 
(a)_{\cal N}\gamma (b)_{\cal N}$ also holds.  \\
The Theorem 2.6 in [1] can be read as follows:\smallskip

\noindent{\bf Proposition.} {\it Let M be a $po$-$\Gamma$-semigroup 
and $a,b\in M$. Then we have the following:

$(1)$ If $(a)_{\cal N}\preceq (b)_{\cal N}$, then $N(b)\subseteq 
N(a)$.

$(2)$ If $N(b)=N(a)$, then $(a)_{\cal N}\preceq (b)_{\cal 
N}$.}\medskip

\noindent{\bf Proof.} The proof of (2): Let $N(b)=N(a)$ and 
$\gamma\in \Gamma$. Since $(a,b)\in\cal N$, we have $(a\gamma 
b,b\gamma b)\in\cal N$. Since $(b\gamma b,b)\in\cal N$, we get 
$(a\gamma b, b)\in\cal N$, then $(a\gamma b)_{\cal N}=(b)_{\cal 
N}=(a)_{\cal N}$, and $(a)_{\cal N}\preceq (b)_{\cal N}$. 
$\hfill\Box$\\
This is Theorem 2.7 in [1]:\smallskip

\noindent{\bf Theorem 2.7.} {\it Let $a\in M$. Then the following 
sets are equal:

$(1)$ $K(a)=\{b\in M \mid (b)_{\cal N} \succ (a)_{\cal N}\}$.

$(2)$ $A=\{b\in N(a) \mid (b)_{\cal N} \succ (a)_{\cal N}\}$.

$(3)$ $B=\bigcup {\{(b)_{\cal N} \mid (b)_{\cal N} \succ (a)_{\cal 
N}}\}$.

$(4)$ $C=N(a)\backslash (a)_{\cal N}$.}

To prove that $K(a)\subseteq A$ or $B\subseteq C$ is not necessary to 
say that $N(b)\varsubsetneq N(a)$ to have $b\in N(a)$. In the proof 
of $C\subseteq A$ the ``$\Leftarrow$'' part of Theorem 2.6 has been 
used. So we do not know if the Theorem 2.7 in [1] is true. The 
Corollary 2.8 of the paper is based on Theorem 2.7, so we cannot say 
that that corollary is true as well. The Corollary 2.8 has been also 
used in Corollary 2.10 and in the Example 2.9.

Before the Theorem 2.11, the author wrote: ``To prove the following 
theorems we use some important notions and results proved in [7, 
1.3.2] for ordered semigroups, the modification of which can be 
easily done for the ordered $\Gamma$-semigroups". The [7, 1.3.2] 
mentioned above is the Reference [2] of the present note but there 
are no such results in the book by J. Howie. In the Proposition 1.3.2 
Howie showed the well known that every lower semilattice is an 
idempotent and commutative semigroup and conversely, that every 
commutative  and idempotent semigroup is a lower semilattice (which 
is true for upper semilattices as well).\\
And this is Theorem 2.11 in [1]: \medskip

\noindent{\bf Theorem 2.11.} {\it The following statements are 
equivalent:

$(1)$ M is a semilattice such that $a\le a\gamma a$ for every $a\in 
M$ and every $\gamma\in\Gamma$.

$(2)$ For every $a\in M$, $N(a)=[a)$.

$(3)$ $\cal N$ is the equality relation on M}.\smallskip

First of all a $\Gamma$-semigroup $M$ is a semilattice means that the 
$\Gamma$-semigroup has the properties $a\gamma a=a$ and $a\gamma 
b=b\gamma a$ for every $a,b\in M$ and every $\gamma\in\Gamma$. So 
property (1) is not true as there is no idempotent (and so idempodent 
and commutative) $\Gamma$-semigroup such that $a< a\gamma a$ for 
every $a\in M$ and every $\gamma\in\Gamma$. The same is mentioned in 
the introduction of the paper as well, where the author wrote: 
``Also, we will consider a structure of principal filter on ordered 
$\Gamma$-semigroups and by using the relation $\cal N$, we will 
observe that $\cal N$ on any ordered $\Gamma$-semigroup $M$ is the 
equality relation on $M$ if and only if $M$ is a semilattice having 
the property $a\le a\gamma a$ for all $a\in M$, $\gamma\in\Gamma$". 
However, in the proof of $(3)\Rightarrow (1)$ of the same theorem he 
shows that $a\gamma a=a$ and $a\gamma b=b\gamma a$ for every $a,b\in 
M$ and  says ``this shows that $M$ is a semilattice as required".

In the proof of the implication $(1)\Rightarrow (2)$, to prove that 
$[a)$ is a filter containing $a$, the author considers $b,c\in M$ 
such that $b\gamma c\in [a)$ ``for all $\gamma\in\Gamma$" (and uses 
the ``for all $\gamma\in\Gamma$" in the rest of the proof as well). 
The correct is: Suppose $b,c\in M$ and $\gamma\in \Gamma$ such that 
$b\gamma c\in [a)$.
From the fact that $M$ is a semilattice, he concludes that
there exist $\gamma_1,\gamma_2\in \Gamma$ such that $b=b\gamma _1 b$ 
and $c=c\gamma_2 c$ which is wrong. Then he wrote: ``Since $b\gamma 
c\ge a$ for all $\gamma\in \Gamma$ (for which we already said is not 
true), we have $b\gamma_1 c\ge a$ and there exists $\gamma_2\in 
\Gamma$ such that $a=a\gamma_2b\gamma_1 c$'' (which is also wrong as 
the order on $M$ does not have this property). This being wrong, the 
rest of the proof that
$$a\gamma_1 b=b\gamma_1 a=a\gamma_2 b\gamma_1 c\gamma_1 b=
a\gamma_2 b\gamma_1 b\gamma_1 c=a$$and that
there exists $\gamma_3\in \Gamma$ such that $a=a\gamma_3 b\gamma_2 
c$, hence
$$a\gamma_2 c=c\gamma_2 a=a\gamma_3b\gamma_2 c\gamma_2 
c=a\gamma_3b\gamma_2c=a,$$cannot be true. Besides, from $a\gamma_1 
b=a$ and from $a\gamma_2 c=a$ we cannot conclude that $b\ge a$ and 
$c\ge a$ to have $b,c\in [a)$.

In the proof of the implication $(3)\Rightarrow (1)$, the phrase 
``since $(a)_{\cal N}$ and $(b)_{\cal N}$ are both semilattice 
congruence classes on $M$" should be replaced by ``since ${\cal N}$ 
is a semilattice congruence on $M$". The $(a)_{\cal N}\gamma 
(a)_{\cal N}\subseteq (a)_{\cal N}$ should be replaced by
$(a)_{\cal N}\gamma (a)_{\cal N}=(a)_{\cal N}$. Since $(a)_{\cal 
N}\gamma (a)_{\cal N}=(a)_{\cal N}$ and $(a)_{\cal N}\gamma (b)_{
\cal N}=(b)_{\cal N}\gamma (b)_{\cal N}$ we have $a\gamma a=a$ and 
$a\gamma b=b\gamma a$ and the author says ``this shows that $M$ is a 
semilattice as required", which is right but in contrast to the 
property (1) of the theorem where the ``$a\le a\gamma a$ $\forall\; 
a\in M$, $\gamma\in\Gamma$" was added. From the proof of the theorem 
it is clear that the author has misunderstood the Proposition 1.3.2 
in [2] the contain of which is mentioned above.\medskip

\noindent We denote by $[a)$ the subset of $M$ defined by 
$[a):=\{t\in M \mid t\ge a\}$.\\
A $po$-$\Gamma$-groupoid $M$ is said to be a {\it band} if $a\gamma 
a=a$ for every $a\in M$ and all $\gamma\in\Gamma$.\medskip

\noindent{\bf Proposition.} {\it Let $M$ be a $po$-$\Gamma$-groupoid 
in which the order ``$\le$" has the following property: $$a\le b 
\;\Longrightarrow\; a\gamma b=a \;\, \forall \; \gamma\in 
\Gamma.$$Then $M$ is a band}.\smallskip

 \noindent{\bf Proof.} Let $a\in M$. Take an element 
$\gamma\in\Gamma$ $(\Gamma\not=\emptyset)$. Since $a\le a$, by 
hypothesis, we have $a\gamma a=a$.$\hfill\Box$\\
\noindent This is the corrected form of the Theorem 2.11 in 
[1]:\medskip

\noindent{\bf Theorem 2.11.} {\it Let M be a $po$-$\Gamma$-semigroup. 
Then we have the following:
\begin{enumerate}
\item[$(1)$] If M is a band then, for any $a\in M$, the set $[a)$ 
is a subsemigroup of M.
\item[$(2)$] If $N(a)=[a)$ for every $a\in M$, then the relation 
$\cal N$ is the equality relation on M.
\item[$(3)$] If $\cal N$ is the equality relation on M, then M is 
a semilattice.
\item[$(4)$] In particular, if M is a commutative band (i.e. a 
semilattice) and the order ``$\le$" on M satisfies the 
relation $$a\le b \;\Longleftrightarrow\; a\gamma b=a \;\, 
\forall \; \gamma\in \Gamma$$then, for every $a\in M$, we 
have $N(a)=[a)$.
\end{enumerate} }\smallskip

\noindent{\bf Proof.} (1) First of all, $[a)$ is a nonempty subset of 
$M$ as $a\in [a)$. Let $x,y\in [a)$ and $\gamma\in\Gamma$. Since 
$x\ge a$, $y\ge a$, $\gamma\in\Gamma$ and $\Gamma$ is a band, we have 
$x\gamma y\ge a\gamma a=a$, so $x\gamma y\in [a)$. Thus $[a)$ is a 
subsemigroup of $M$.

(2) Let $(a,b)\in\cal N$. Since $a\in N(a)=N(b)=[b)$, we have $a\le 
b$. Since $b\in N(b)=N(a)=[a)$, we have $b\le a$. Thus we have $a=b$, 
and ${\cal N}$ is the equality relation on $M$.

(3) Let $a,b\in M$ and $\gamma\in\Gamma$. Since the relation $\cal N$ 
is a semilattice congruence on $M$, we have $(a\gamma a,a)\in\cal N$ 
and $(a\gamma b,b\gamma a)\in\cal N$. Since $\cal N$ is the equality 
relation on $M$, we have $a\gamma a=a$ and $a\gamma b=b\gamma a$, 
thus $M$ is a semilattice.

(4) Let $a\in M$. As we have already seen in (1), the set $[a)$ is a 
subsemigroup of $M$. Let $x,y\in [a)$ and $\gamma\in\Gamma$ such that 
$x\gamma y\in [a)$. Then $x\in [a)$ and $y\in [a)$. Indeed: Let 
$\mu\in\Gamma$ $(\Gamma\not=\emptyset)$. Since $x\gamma y\ge a$, by 
hypothesis, we have $a\mu(x\gamma y)=a$. Since $M$ is a band, we 
have\begin{eqnarray*}a\mu x&=&{\Big(}a\mu (x\gamma y){\Big )}\mu x
=a\mu {\Big(}(x\gamma y)\mu x{\Big )}=a\mu {\Big(}x\mu (x\gamma 
y){\Big )}=(a\mu x)\mu (x\gamma y)\\&=&a\mu (x\mu x)\gamma y=a\mu 
x\gamma y=a,\end{eqnarray*}and$$a\mu y=(a\mu x\gamma y)\mu y=a\mu 
x\gamma (y\mu y)=a\mu x\gamma y=a,$$ so $x\ge a$ and $y\ge a$, that 
is $x,y\in[a)$. If $x\in [a)$ and $M\ni y\ge x$, then $y\in [a)$. 
Indeed: Since $x\in [a)$, we have $x\ge a$. Then $y\ge x\ge a$, so 
$y\in [a)$. Let $T$ be a filter of $M$ such that $a\in T$. Then 
$[a)\subseteq T$. Indeed: If $x\in [a)$, then $x\in M$ and $x\ge a\in 
T$. Since $T$ is a filter of $M$, we have $x\in T$. As $[a)$ is the 
smallest (with respect to the inclusion relation) filter of $M$ 
containing $a$, we have $N(a)=[a)$.$\hfill\Box$

Taking into account the Proposition above we notice that if $M$ is a 
commutative $po$-$\Gamma$-semigroup satisfying the relation $$a\le b 
\;\Longleftrightarrow\; a\gamma b=a \mbox { for all } \gamma 
\in\Gamma,$$ then $M$ is a commutative band (and so a semilattice) 
and, for every $a\in M$, we have $N(a)=[a)$.\\
This is the Theorem 2.13 in [1]:\medskip

\noindent{\bf Theorem 2.13.} {\it Let $\sigma$ be a complete 
semilattice congruence on an ordered $\Gamma$-semigroup $M$ and $Y$ 
the semilattice $M/{\sigma}$. Then for any $x\in Y$, we have
\begin{enumerate}
\item[$(1)$] $M_x$ is the union of some ${\cal N}$-classes.

\item[$(2)$] The set $T=\bigcup \{M_y \mid y\succeq x, y\in Y\}$ 
is a filter.

\item[$(3)$] For any $a\in M_x$, $N(a)=T$ if and only if $\sigma$ 
is the smallest complete semilattice congruence on $M$.
\end{enumerate}}
The relation $``y\succeq x"$ for the elements of $M/\sigma$ has not 
defined in the paper. On p. 194 in [1] the definition is only for 
$\sigma={\cal N}$.

We define $(a)_\sigma \preceq (b)_\sigma \Longleftrightarrow 
(a)_\sigma=(a)_\sigma\gamma (b)_\sigma:=(a\gamma b)_\sigma$ for all 
$\gamma\in\Gamma$.

The $M_x$ has not defined in the paper, apparently it is the 
$(x)_\sigma$. As far as the property (2) is concerned, it should be 
clarified in (2) if the filter mentioned in it was in $M$ or in 
$M/{\sigma}\;$. It seems that it is in $M/{\sigma}$ while according 
to the proof of (2) it is in $M$. In fact, the author tried to prove 
in (2) that ``$T$ is a subsemigroup of $M$" and that if ``$a,b\in M$ 
and $\gamma\in\Gamma$ such that $a\gamma a\in T$, then $a\in T$ and 
$b\in T$", which means that he considers the filter in $M$. The proof 
contains serious mistakes in it. For example, in several parts of the 
proof $M_z\Gamma M_t\subseteq M_{z\Gamma t}$ is written. \\
The Theorem 2.13 in [1] could be replaced by the following:\medskip

\noindent{\bf Theorem 2.13.} {\it Let M be a $po$-$\Gamma$-semigroup 
and $\sigma$ a semilattice congruence on M. Then the following 
property is satisfied:

If $x\in M$, then the set $T:=\{(y)_\sigma \mid y\in M, 
(y)_\sigma\succeq (x)_\sigma\}$ is a filter in $M/\sigma$.}\medskip

\noindent{\bf Proof.} Take an element $\gamma\in\Gamma$ 
$(\Gamma\not=\emptyset)$. Since $\sigma$ is a semilattice congruence 
on $M$, we have $(x,x\gamma x)\in\sigma$, then $(x)_\sigma=(x\gamma 
x)_\sigma$, so $(x)_\sigma\succeq (x)_\sigma$, and $(x)_\sigma\in T$. 
Thus $T$ is a nonempty subset of $M/\sigma$. Let now $(y)_\sigma, 
(z)_\sigma\in T$ and $\gamma\in\Gamma$. Then 
$(y)_\sigma\gamma(z)_\sigma\in T$, that is $(y\gamma z)_\sigma\in T$. 
Indeed: Since $(y)_\sigma\succeq (x)_\sigma$, $(z)_\sigma\succeq 
(x)_\sigma$ and $\sigma$ is a semilattice congruence on $M$, we 
have$$(y\gamma z)_\sigma=(y)_\sigma\gamma(z)\succeq (x)_\sigma\gamma 
(x)_\sigma=(x\gamma x)_\sigma=(x)_\sigma,$$so $(y\gamma z)_\sigma\in 
T$. \noindent Let $a,b\in M$ and $\gamma\in\Gamma$ such that 
$(a)_\sigma \gamma (b)_\sigma\in T$. We have to prove that 
$(a)_\sigma\in T$ and $(b)_\sigma\in T$, that is $(a)_\sigma\succeq 
(x)_\sigma$ and $(b)_\sigma\succeq (x)_\sigma$, which means that 
$(x)_\sigma=(x\mu a)_\sigma$ and $(x)_\sigma=(x\mu b)_\sigma$ for 
every $\mu\in\Gamma$. \\Let now $\mu\in\Gamma$. By hypothesis, we 
have $(a\gamma b)_\sigma\in T$, then $(a\gamma b)_\sigma\succeq 
(x)_\sigma$, and so $(x)_\sigma=(x\xi a\gamma b)_\sigma$ for every 
$\xi\in\Gamma$. Thus we have $(x)_\sigma=(x\mu a\gamma b)_\sigma$. 
Then we get$$(x\mu a)_\sigma=(x)_\sigma\mu (a)_\sigma=(x\mu a\gamma 
b)_\sigma\mu (a)_\sigma=(x\mu a\gamma b\mu a)_\sigma.$$Since $\sigma$ 
is a semilattice congruence on $M$, we have ${\Big(}(a\gamma b)\mu a, 
a\mu (a\gamma b){\Big)}\in\sigma$ and $(a\mu a,a)\in\sigma$. Then 
${\Big(}(a\mu a)\gamma b, a\gamma b{\Big)}\in\sigma$, and $(a\gamma 
b\mu a, a\gamma b)\in\sigma$. Since $(a\gamma b\mu a, a\gamma 
b)\in\sigma$ and $\sigma$ is a congruence on $M$, we have $(x\mu 
a\gamma b\mu a, x\mu a\gamma b)\in\sigma.$ Hence we obtain $$(x\mu 
a)_\sigma=(x\mu a\gamma b)_\sigma=(x)_\sigma.$$ Similarly we prove 
that $(x\mu b)_\sigma=(x)_\sigma$ for all $\mu\in\Gamma$. Finally, 
let $(y)_\sigma\in T$ and $M/\sigma\ni (z)_\sigma\succeq (y)_\sigma$.  
Then $(z)_\sigma\in T$. Indeed: Take an element $\gamma\in\Gamma$ 
$(\Gamma\not=\emptyset)$. Since $(y)_\sigma\succeq (x)_\sigma$ and
$(z)_\sigma\succeq (y)_\sigma$, we have $(x)_\sigma=(x)_\sigma\gamma 
(y)_\sigma$ and $(y)_\sigma=(y)_\sigma\gamma (z)_\sigma$. Hence we 
obtain $$(x)_\sigma=(x)_\sigma\gamma (y)_\sigma=(x)_\sigma\gamma 
{\Big(}(y)_\sigma\gamma (z)_\sigma{\Big)}={\Big(}(x)_\sigma\gamma 
(y)_\sigma{\Big)}\gamma (z)_\sigma=(x)_\sigma\gamma (z)_\sigma,$$so 
$(z)_\sigma\succeq (x)_\sigma$, and $(z)_\sigma\in 
T$.$\hfill\Box$\medskip

The Corollaries 2.14 and 2.15 of the paper are based on Theorem 2.13. 
The Example 2.16 is based on Theorem 2.11. In Examples 1.3 and 1.6 
the author defines an order on ``$\Gamma$" while for a 
$\Gamma$-semigroup, the set $\Gamma$ is just a nonempty set and not 
an ordered set.

Finally it should be noted that except of the case in which we search 
for a counterexample in which case an example for an ordered 
semigroup is enough, examples of ordered $\Gamma$-semigroups in which 
the set $\Gamma$ consists only by one element are actually examples 
of ordered semigroups. A sufficient example of an ordered 
$\Gamma$-semigroup should be an example in which the set $\Gamma$ has 
at least two elements. The examples of the paper in [1] are examples 
of ordered semigroups. A counterexample is given in the Example 1.11, 
but this also being an example of an ordered semigroup, just looking 
at the table and the figure of it one can immediately concludes that 
this is an example of a commutative ordered semigroup (and so of a 
commutative ordered $\Gamma$-semigroup as well) while the author gets 
the assumption that it is not a $\Gamma$-semigroup to prove that it 
is.

 Many of the results of the paper in [1] hold in 
$po$-$\Gamma$-groupoids in general.

\end{document}